\documentclass[pdflatex,10pt]{article}
\usepackage[utf8]{inputenc}

\usepackage{amsmath,color}
\usepackage{amssymb,amsthm}
\usepackage{subfig}
\usepackage[permil]{overpic}
\usepackage[multiple]{footmisc}
\usepackage{xr}
\usepackage[left=1.6cm,right=1.6cm,top=2.50cm,bottom=2.50cm]{geometry}
\usepackage{graphicx}
\usepackage[font=small,labelfont=bf,
   justification=justified,
   format=plain,labelsep=space]{caption}
\usepackage{subcaption}
\usepackage{float}
\usepackage[english]{babel}
\usepackage[font=small,labelfont=bf,justification=centering]{caption}
\usepackage[T1]{fontenc}
\usepackage{lastpage}
\usepackage{enumerate}
\usepackage{enumitem}
\usepackage{lmodern}
\usepackage{array}
\usepackage{bm}
\usepackage{multicol}
\usepackage{multirow}
\usepackage{dsfont}
\usepackage{tensor}
\usepackage{fancyhdr}
\usepackage{listings}
\usepackage{dsfont}
\usepackage{siunitx}
\usepackage{titling}
\usepackage{lipsum}
\usepackage{tabularx}
\usepackage{verbatim}
\usepackage{authblk}
\usepackage{csquotes}
\usepackage{bbm} %for indicator function
\usepackage{hyperref} % good for references
\usepackage{ctable} % for \specialrule command for tables
\usepackage{multirow} % Merge rows in tables 
\usepackage{titlesec} % Modify font sections, subsections, etc.
\usepackage{marvosym}
\usepackage{relsize,exscale} % for big integrals
\usepackage{setspace}
\usepackage[
backend=biber,
dashed=false,
citestyle=numeric,
bibstyle=apa,
giveninits=true,
sorting=anyt
]{biblatex}

\makeatletter
\input{numeric.bbx}
\makeatother

\allowdisplaybreaks

\graphicspath{{Fig/}}

\addto\captionsenglish{}

\titleformat{\subsection}{\normalfont\large\raggedright\it}{\thesubsection}{1em}{}

\makeatletter
\newcommand{\thickhline}{%
    \noalign {\ifnum 0=`}\fi \hrule height 1pt
    \futurelet \reserved@a \@xhline
}
\newcolumntype{"}{@{\hskip\tabcolsep\vrule width 1pt\hskip\tabcolsep}}
\makeatother

\usepackage{adjustbox}

\makeatletter
\newsavebox{\@brx}
\newcommand{\llangle}[1][]{\savebox{\@brx}{\(\m@th{#1\langle}\)}%
  \mathopen{\copy\@brx\kern-0.5\wd\@brx\usebox{\@brx}}}
\newcommand{\rrangle}[1][]{\savebox{\@brx}{\(\m@th{#1\rangle}\)}%
  \mathclose{\copy\@brx\kern-0.5\wd\@brx\usebox{\@brx}}}
\makeatother

\theoremstyle{definition}

\numberwithin{equation}{section}

\usepackage{graphicx} % Required for inserting images

\newcommand{\eps}{\varepsilon}

\title{Large Deviation Minimisers for Stochastic Partial Differential Equations with Degenerate Noise}
\author[1,\dag]{Paolo Bernuzzi}
\author[2]{Tobias Grafke}
\affil[1]{\footnotesize{Technical University of Munich, School of Computation Information and Technology, Department of
Mathematics, Boltzmannstraße 3, 85748 Garching, Germany}}
\affil[2]{\footnotesize{Mathematics Institute, University of Warwick, Coventry CV4 7AL, United Kingdom}}
\affil[$\dag$]{Author to whom any correspondence should be addressed.\smallskip 

Email addresses: paolo.bernuzzi@ma.tum.de (P. Bernuzzi), T.Grafke@warwick.ac.uk (T. Grafke).}
\date{\today}

\begin{document}

\maketitle

% REQUIRED

\begin{abstract}
  Noise-induced transitions between multistable states happen in a
  multitude of systems, such as species extinction in biology, protein
  folding, or tipping points in climate science. Large deviation
  theory is the rigorous language to describe such transitions for
  non-equilibrium systems in the small noise limit. At its core, it
  requires the computation of the most likely transition pathway,
  solution to a PDE constrained optimisation problem. Standard methods
  struggle to compute the minimiser in the particular coexistence of (1)
  multistability, i.e.~coexistence of multiple long-lived states, and (2)
  degenerate noise, i.e.~stochastic forcing acting only on a small
  subset of the system's degrees of freedom. In this paper, we
  demonstrate how to adapt existing methods to compute the large
  deviation minimiser in this setting by combining ideas from optimal
  control, large deviation theory, and numerical optimisation. We show
  the efficiency of the introduced method in various applications in
  biology, medicine, and fluid dynamics, including the transition to
  turbulence in subcritical pipe flow.
\end{abstract}


\small{\textbf{MSC codes:} 60H15,   % Stochastic partial differential equations
65C50,   % Other Computational Problems in Prob
65Z05,   % Application to Physics (Numerical Analysis
82B26.    % Phase Transitions (general) in equilibrium stat mech
%,82B31  % Stochastic methods (Eq Stat Mech)}

\small{\textbf{Funding:} This work was supported by the European Union’s Horizon 2020 research and innovation programme under Grant Agreement 956170. TG acknowledges support from EPSRC projects
\mbox{EP/T011866/1} and \mbox{EP/V013319/1}.}}

% Header
\pagestyle{fancy}
\fancyhead{}
\renewcommand{\headrulewidth}{0pt}
\fancyhead[C]{\textit{Large Deviation Minimisers for Degenerate SPDEs}}

\section{Introduction}

Rare events are occurrences that, by nature, happen
infrequently. Examples can be found in biology, chemistry
\cite{rubino2009rare}, climate science~\cite{crommelin2015stochastic}
and material science~\cite{liu2015efficient}. Despite their rarity,
their impact is often qualitatively significant and far-reaching. For
instance, the prediction of extreme heat waves
\cite{ragone2021rare,ragone2018computation} holds importance due to
the economic impact and health risks that the events pose. Moreover, the forecast of rare events in climate science is connected with the ever-growing topic of tipping points \cite{lenton2013origin}, and the fact that many anomalies are more likely to happen due to human intervention. In recent years, the causes of different sensible phenomena have been deeply studied along with the possible consequences of their tipping \cite{loriani2023tipping}. Their severe global effects and the expectation of tipping cascades underlines the importance of the construction of valid early-warning signs \cite{bernuzzi2024warning,bernuzzi2023bifurcations,dakos2024tipping}.

The quantitative treatment of rare events requires the introduction of
a stochastic model of the system at hand. In such models, rare but
sudden state changes are typically characterised by multistability,
i.e.~the coexistence of multiple long-lived states, and the
catastrophic rare event corresponds to a noise-induced transition
between two such states. The probability of qualitative state changes
is related to the stability of a steady state
\cite{blomker2005modulation} and to the intensity of the noise~\cite{berglund2023stochastic,bernuzzi2023bifurcations}. Rare event 
methods based on sampling rely on biasing techniques to increase the 
likelyhood of a rare event, such as importance splitting 
techniques~\cite{cerou-guyader:2007, lestang-ragone-brehier-etal:2018, cerou-guyader-rousset:2019}.
In this paper, we consider the noise strength to be infinitesimally
small. This puts us in the purview of Freidlin-Wentzell theory of
large deviations~\cite{freidlin1998random}, meaning that both the
limiting event probability scaling and its most likely pathway of
occurence are rigorously predictable, and available in a sampling-free
and determinstic manner. Obtaining the most likely pathway provides insights
into the physical mechanism of the transition as it crosses between
two different deterministic basins of attraction. The analytic
properties of such a solution have been described for certain
stochastic differential equations (SDEs)
\cite{breden2019rigorous,xu2014rare}, and the suitable technique
involved in its numerical description usually depends on the
particular model
\cite{bisewski2019rare,grafke2019numerical,khoo2019solving}. Especially
challenging is the assumption of non-white noise, which can be found
in many applications, such as privacy risk management
\cite{borovykh2023privacy}. In this paper, we focus on the study of
white-in-time Gaussian but \emph{degenerate} noise, which forces only
specific components of the system. This type of noise is fairly
generically assumed in many scientific situations, where stochasticity
is inserted to model a subset of unknown or unknowable degrees of
freedom~\cite{blomker2012numerical}. For example, one would force
only small or large-scale structures in fluid
turbulence~\cite{dong2018ergodicity}, or only surface freshwater
influx or wind stresses in ocean
modelling~\cite{baars2017continuation,bernuzzi2024warning}. As we will
see, the coexistence between multistability and degenerate noise poses
a challenge to numerical algorithms to compute the most likely
transition trajectory, and even more so if it is multiplicative at the
same time.

The paper is structured as follows. In section~\ref{sec:FW theory},
the Freidlin-Wentzell theory is introduced and the most likely path
that describes a rare event is obtained as a minimiser of the action
functional on a manifold. Subsequently, in section~\ref{sec:AS
  method}, we introduce the adjoint state method, a standard method in
a numerical optimisation problem~\cite{plessix2006review}, as well as
its application to the construction of the path. This includes a
comparison with other techniques, and a discussion of its ability to
handle degenerate noise. Lastly, a selection of applications are shown
in section~\ref{sec:APP}, with particular focus on (degenerate noise,
multistable) stochastic partial differential equations~(SPDEs). Among
the phenomena treated are spike merging in the Gierer-Meinhardt
model~\cite{wei2013mathematical}, pulse initiation in the
FitzHugh-Nagumo model~\cite{beyn2018computation}, and turbulent puff
splitting event for subcritical pipe flow in the Barkley
model \cite{barkley2016theoretical}. The different types of noise
treated are distinguished into additive and multiplicative degenerate
noise. The assumption of boundary noise is also studied.

\section{Freidlin-Wentzell theory}
\label{sec:FW theory}

In this section, we introduce Freidlin-Wentzell theory~\cite{freidlin1998random} and its application to SPDEs. Labelling the domain $\mathcal{X}\subset \mathbb{R}^N$, and with $\eps>0$, we define the cylindrical Wiener process $W$, such that 
\begin{align*}
    W(x,t) = \sum_{i=1}^\infty b_i(x) \xi_i(t)\;,
\end{align*}
for $\left\{b_i\right\}_{i\in\mathbb{N}_{>0}}$ a basis of $L^2(\mathcal{X})$, $\left\{\xi_i\right\}_{i\in\mathbb{N}_{>0}}$ a collection of Brownian processes, and for any $x\in\mathcal{X}$ and $t\geq0$.
For initial condition $u_0\in L^2(\mathcal{X})$, deterministic drift term $b:L^2(\mathcal{X})\to L^2(\mathcal{X})$ and noise diffusion operator $\sigma:L^2(\mathcal{X})\to L^2(\mathcal{X})$, we assume \cite{da2014stochastic} that there exists a solution $u^\eps$ of the system
\begin{align*}
    \left\{\begin{alignedat}{2}
            &\text{d}u^\eps(x,t) = b\left(u^\eps(x,t)\right) \text{d}t + \sqrt{\eps} \sigma\left(u^\eps(x,t)\right) \text{d} W(x,t) \quad &&, \; x\in\mathcal{X}\;, \; t\geq 0 \;,\\
            &u^\eps(x,0) = u_0(x) &&, \; x\in\mathcal{X}\;.
    \end{alignedat}\right.
\end{align*}
For a fixed time $T>0$, we introduce the action functional
\begin{align*}
    I_T(\phi) = \underset{\left\{g\in H^1: \phi(t)=u_0+\int_0^t b(\phi(s)) \text{d}s + \int_0^t \sigma(\phi(s)) \dot{g}(s) \text{d}s, t\leq T \right\}}{\text{inf}} \frac{1}{2} \int_0^T \lvert\lvert \dot{g}(s) \rvert\rvert^2 \text{d}s,
\end{align*}
which quantifies the effect of the noise on the solution path $\phi$ in the time interval $[0,T]$. Then~\cite{freidlin1998random}, under the assumption of small noise intensity $0<\eps\ll1$, the following exponential estimate holds
\begin{align*}
    \text{log} \left(\mathbb{P} \left( \underset{t \in [0,T]}{\text{sup}} \lvert\lvert u^\eps(t) - \phi(t) \rvert\rvert^2 < \delta \right) \right) \propto - \eps^{-1} I_T(\phi)\;,
\end{align*}
for small enough $\delta>0$. In the limit $\eps\to0$, the most
likely paths $\phi$ minimise the action functional. For rare
transitions between metastable
states~\cite{grafke2019numerical,rolland2016computing}, particular
importance is given to paths in
\begin{align*}
    \mathcal{C}:=\left\{\phi\in C([0,T]) \;|\; \phi(0) = u_0, \phi(T) = u_T \right\},
\end{align*}
for $u_0, u_T\in L^2(\mathcal{X})$ in disjoint deterministically
invariant subsets of $\mathcal{X}$. Under the assumption of small
noise intensity, paths that cross the two states are unlikely by
construction. For $\eps\to 0$ and fixed states $u_0,u_T\in
L^2(\mathcal{X})$, the \emph{instanton} $\bar{\phi}$ is the most likely
solution $u^\eps$ such that $\phi\in\mathcal{C}$, and equivalently, it
satisfies
\begin{align}
  \label{eq:FW-optim}
    \bar{\phi} = \underset{\phi\in\mathcal{C}}{\text{argmin}} \; I_T(\phi) \;.
\end{align}
The nature of the noise is defined by the noise diffusion operator
$\sigma$ and the noise covariance $a:=\sigma^\ast \sigma$, with
${}^\ast$ indicating the adjoint operator in the Hilbert space
$L^2(\mathcal{X})$. If the noise is independent of the space variable
$u^\eps\in L^2(\mathcal{X})$ we call the noise \emph{additive}, while in the
case where $\sigma$ depends on $u^\eps$, the noise is
\emph{multiplicative}. Further, if there exist $u,v\in
L^2(\mathcal{X})$ and $v\not\equiv 0$ , such that $\sigma(u)v=0$, then the
noise is \emph{degenerate}. For all these cases, we define the
Hamiltonian functional
\begin{align*}
    H(\phi,\theta) = \langle b(\phi), \theta \rangle + \frac{1}{2} \langle \theta, a(\phi) \theta \rangle
\end{align*}
for $\theta$ the conjugate momentum of $\phi$. The instanton solves~\cite{grafke2019numerical} the corresponding Hamilton equations
\begin{align*}
    \left\{ \begin{alignedat}{3}
        \dot{\phi} &= \partial_\theta H(\phi,\theta) &&= b(\phi) + a(\phi) \theta &&, \\
        \dot{\theta} &= -\partial_\phi H(\phi,\theta) &&= -\partial_\phi \left(b(\phi)\right) \theta - \frac{1}{2} \theta \;\partial_\phi \left(a(\phi)\right) \theta &&.
    \end{alignedat}
    \right.
\end{align*}
The first Hamilton equation indicates that the weighted conjugate
variable $\sigma(\phi) \theta$ can be interpreted as the optimal noise
on the instanton path. Yet, in the case of degenerate noise, it is
clear that its role is more subtle. In particular, we know that in
this case there must be modes $v\in L^2(\mathcal{X})$ that remain unforced
regardless of the choice of $\theta$. In other words, while it remains
true that $\sigma(\phi)\theta$ is the optimal noise, the mapping
between noise and trajectory $\phi$ is no longer one-to-one: there
exist trajectories $\phi$ that cannot be realised by any value of the
noise.

\medskip
This fact stands in the way of applying traditional methods for
computing large deviation minimisers, such as the minimum action
method (MAM)~\cite{e-ren-vanden-eijnden:2004, zhou-ren-e:2008} and its
geometric counterpart~\cite{grafke-schaefer-vanden-eijnden:2017,
heymann-vanden-eijnden:2008-a,
vanden-eijnden-heymann:2008}. Such methods obtain the
optimal transition path of the Lagrangian optimisation problem via
relaxation or gradient descent, which necessitates inverting the
noise covariance matrix $a(\phi)$. In other words, methods that rely
on optimising the path instead of the noise are inadequate in the
presence of degenerate noise and non-invertible covariance
matrices.

On the other hand, methods based on the Hamiltonian formalism, as
introduced above, do not suffer from the same shortcoming, as they
circumvent inverting the noise covariance by considering not the
Lagrangian optimisation problem, but its Legendre transform in the
form of the Hamilton's equations. Corresponding methods have been
employed successfully in the presence of degenerate noise even for
high-dimensional systems~\cite{grafke-grauer-schaefer-etal:2014,
  grafke2019numerical, simonnet:2023}. Unfortunately, with the notable exception 
of~\cite{zakine-vanden-eijnden:2023}, these approaches are in general
incapable of dealing with metastability, which always implies
nonconvexity of the rate function, leading to nonuniqueness of the
boundary conditions of the adjoint
variable~\cite{alqahtani-grafke:2021}.

In this paper, we suggest to unify these two approaches, by
simultaneously applying the adjoint state method (allowing us to
consider degenerate noise), while also convexifying the problem
through the application of the augmented Lagrangian
method~\cite{nocedal1999numerical} or penalisation of the endpoint
constraint (see also the equivalent idea of the generalised canonical
ensemble in statistical
mechanics~\cite{costeniuc-ellis-touchette-etal:2005}). This
combination allows us to deal with metastability and degenerate noise
simultaneously, and to do so for high-dimensional and complex systems
such as stochastic PDEs.

\section{Adjoint state method}
\label{sec:AS method}

The adjoint state method is a well-known technique for the solution of
constained minimisation problems~\cite{plessix2006review}. Its
application in the action minimisation problem is easily justifiable
and analyzable when realizing that the Freidlin-Wentzell minimisation
problem~\eqref{eq:FW-optim} can be interpreted as minimising the
$L^2(a, \mathcal X)$-cost of the noise under a functional constraint
enforcing all degrees of freedom of the path, both the stochastically
forced and the deterministic ones. Concretely, we define variables
$\phi,\theta$ that assume values in $L^2(\mathcal{X})$ for any
$t\in[0,T]$, and introduce the adjoint state variable, or the Lagrange
multiplier, $\mu$, which also assumes values in $L^2(\mathcal X)$ for
any $t\in[0,T]$. Further, the endpoint constraint is enforced through
a penalty parameter $\lambda>0$. The cost function is then constructed
as follows,
\begin{align}
    \label{eq:adjoint-cost}
    &J(\phi,\theta,\mu,\lambda) =
    \frac{1}{2} \underbrace{\int_0^T \lvert\lvert \sigma(\phi) \theta \rvert\rvert^2 \text{d}t}_{\text{Weighted action}} \\\nonumber
    &+ \underbrace{\int_0^T \left\langle \mu, \dot{\phi}-\partial_\theta H(\phi,\theta) \right\rangle \text{d}t}_{\text{First Hamilton equation enforcer}} + \frac{1}{2} \underbrace{\lambda \lvert\lvert \operatorname{F} \left( \phi(T)-u_T \right)\rvert\rvert^2}_{\text{Penalty term}} \;.
\end{align}
The first term in~\eqref{eq:adjoint-cost} is the actual cost of the
stochastic forcing, appropriately weighed by the noise covariance
$\sigma(\phi)$. The second term employs the adjoint trajectory
$\mu(t)$ to enforce the first Hamilton's equation, which yields the
unique path $\phi(t)$ caused by the noise
$\sigma(\phi)\theta$. Finally, the last term enforces the endpoint
constraint of the path via a penalty term. In other words, the method
corresponds to a Lagrange Multiplier Penalty method
\cite{nocedal1999numerical} in path space. The linear operator
$\operatorname{F}$ is employed as a filter to capture characteristics
of the final condition $u_T$. In the upcoming examples, it is assumed
to be the identity operator $\operatorname{I}$ or a multiplication
operator. The cost function is iteratively minimised, for example via
gradient descent or quasi-Newton method, until the partial derivatives
of $J$ regarding $\phi,\theta$, and $\mu$ have norms below a
pre-chosen threshold. Each iteration is solved in accordance with the
steps to follow:
\begin{enumerate}
    \item We enforce $\partial_\mu J = 0$, i.e.
    \begin{align*}
        \partial_{\mu} J &= \underbrace{\dot{\phi} - \partial_{\theta} H(\phi,\theta)}_{\text{First Hamilton equation}} \\
        &= \dot{\phi} - b(\phi) - \sigma^\ast(\phi) \sigma(\phi) \theta \;,
    \end{align*}
    by solving the first Hamilton equation, for which $\phi=u_0$ is
    chosen as the initial condition. The variable $\theta$ assumes the
    role of the conjugate momentum of $\phi$. For our SPDE-examples,
    this is solved through the mild solution formula
   ~\cite{da2014stochastic}, in the eigenbasis of the Laplacian, while
    for SDEs we employ the implicit Euler-Maruyama method.
    \item
    \small Next, we enforce $\partial_\phi J = 0$, i.e.
    \begin{align*}
        \partial_{\phi} J &= \underbrace{-\dot{\mu} - \partial_{\phi} H(\phi,\mu)}_{\text{Second Hamilton equation}} \\
        &+ \underbrace{\frac{1}{2} \partial_{\phi} \left(\left( \sigma \left(\phi\right) (\theta -\mu) \right)^2\right)}_{\text{Error term}} + \underbrace{\delta_{T} \lambda  \operatorname{F}^* \operatorname{F} \left(\phi-u_T\right)}_{\text{Final condition enforcer}} \;.
    \end{align*}
    Such a constraint is equivalent to solving the second Hamilton
    equation, with $\mu$ in the role of the conjugate variable and
    with an error term. The equation is solved backward in time and
    initialised as $\mu(T)=\lambda \left(\phi(T)-u_T\right)$. Such a
    final condition in the adjoint state variables aims to enforce the
    targeted final condition in the path variable $ F \phi(T) = F
    u_T$. Similarly to the previous step, the numerical computation of
    the SPDEs is resolved through the mild solution formula
   ~\cite{da2014stochastic}, in the eigenbasis of the Laplacian with
    corresponding boundary conditions. The SDEs are solved through the
    implicit Euler-Maruyama method and $\operatorname{F} =
    \operatorname{I}$.
    \item Lastly, the partial derivative
    \begin{align*}
        \partial_{\theta} J = a\left(\phi\right) \left(\theta - \mu \right) = \sigma^\ast(\phi) \sigma(\phi) \left(\theta - \mu \right)
    \end{align*}
    is computed. This gradient is guaranteed to be a descent direction
    of the cost function~\eqref{eq:adjoint-cost} in the space of
    $\theta$-variables, and consequently, it can be used to update the
    variable $\theta$, in our case through the L-BFGS
    method~\cite{nocedal1999numerical}.
\end{enumerate}
The iteration is halted as the norm of $\partial_{\theta} J$ is less
than a fixed tolerance value, $tol>0$. At the actual minimiser, we
have $\partial_{\theta} J\left(\phi,\theta,\mu,\lambda\right)=0$ and it
is implied that $\sigma(\phi) \left( \theta - \mu \right) =
0$. Therefore, the error term in the second step is null, and the
Lagrange multiplier $\mu$ solves the first Hamilton equation in the
first step. This implies that the couple $\left( \phi, \mu \right)$
solves the Hamilton equations and that $\phi$ is the instanton of the
SPDE under initial condition $\phi(0)=u_0$ and final condition
$\operatorname{F} \phi(T)=\operatorname{F} u_T$. Furthermore, this
implies that the optimal noise can be read off as
$\sigma(\phi) \mu(t)$ at convergence.

\medskip
The optimisation problem presented above can be solved for any assumed
transition time $T>0$. Generally, though, the optimisation problem is
harder to solve for longer time intervals. First, this is because
longer time intervals necessitate more numerical timesteps and thus a
larger number of degrees of freedom to optimise. Second, and more
importantly, in a metastable system, the most likely transition will
generally involve a localised-in-time jump from one state to the other
(hence the name 'instanton'). Because of this, for long time
intervals, the optimisation problem is generally very insensitive to
time-translations: the Freidlin-Wentzell action of an earlier or later
jump is almost identical, leading to flat directions in the cost
functional that involve many degrees of freedom and are thus hard to
correct for, even with quasi-Newton methods. On the other hand, this
phenomenologically implies that the optimisation procedure will rather
rapidly converge to the correct transition path, and from then on only
very slowly to the correct transition time.

\section{Applications}
\label{sec:APP}

In this section, we demonstrate the applicability and efficiency of
the above algorithm for the computation of transition pathways in
stochastic complex systems with degenerate noise. In particular, we
will consider systems of increasing complexity, starting with a toy
SDE model of a two-dimensional double-well with degenerate
multiplicative noise in subsection~\ref{sec:two-dimensional-sde},
considering then the Allen-Cahn reaction diffusion SPDE, but forced
only through the boundary in subsection~\ref{sec:one-dimens-allen}, and
the Gierer-Meinhardt SPDE for pattern formation with multiplicative
degenerate noise in subsection~\ref{sec:spike-merging-gierer}. Afterward, we
concentrate concretely on SPDEs of the advection-reaction-diffusion
type involving formation of moving spikes, including the spatially
extended FitzHugh-Nagumo model with additive degenerate noise in subsection~\ref{sec:pulse-init-fitzh}, and lastly the Barkley model for
turbulence proliferation in pipe flow with multiplicative degenerate
noise in subsection~\ref{sec:puff-splitt-barkl}.

\subsection{Two-dimensional SDE with degenerate multiplicative noise}
\label{sec:two-dimensional-sde}

\begin{figure}[h!]
  \centering
  \subfloat{
      \includegraphics[width=0.47\linewidth]{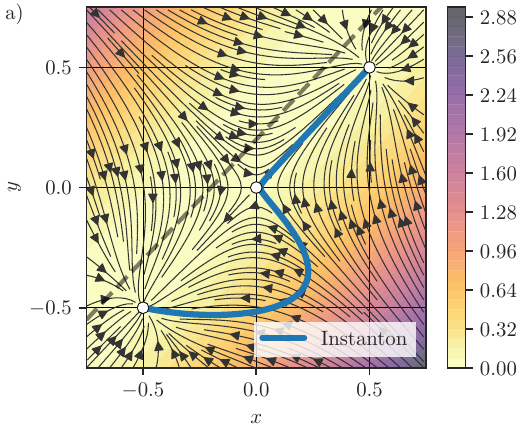}
  }
  \hspace{5mm}
  \subfloat{
      \includegraphics[width=0.47\linewidth]{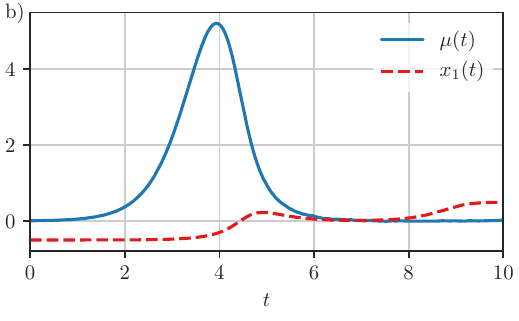}
  }
  \caption{a) Instanton (solid line) for the transition between the
    two stable attractors of the dynamics~\eqref{eq:SDE-2D}
    (streamlines). The shading refers to the strength of the
    multiplicative noise, $a(x_1,x_2)$, with a dashed line at
    $a(x_1,x_2)=0$. The markers indicate the initial point $u_0$, the
    saddle, and the final point $u_T$. b) Lagrange multiplier
    $\mu(t)$ adjoint to the optimal noise (solid blue), and first
    component of the minimiser, $x_1(t)$ (dashed red), with respect to
    time. The figures are obtained for $\lambda=5$ and $tol=10^{-4}$.}
  \label{fig:SDE-2D}
\end{figure}

For $u=(x_1,x_2)\in\mathbb{R}^2$, the SDE
\begin{align}
  \label{eq:SDE-2D}
    \text{d} u^\eps =
    \text{d} \begin{pmatrix}
        x_1 \\ x_2
    \end{pmatrix}
    =
    \begin{pmatrix}
        2 x_2 - (x_1+x_2)^3 \\
        2 x_1 - (x_1+x_2)^3
    \end{pmatrix} \text{d} t + \sqrt{\eps}
    \begin{pmatrix}
         (x_1-x_2+0.2) \\ 0
    \end{pmatrix}
    \text{d} W_t
\end{align}
represents diffusion in a double well potential on the main diagonal
$x_1=x_2$. It has two deterministically stable solutions at
$(x_1,x_2)=u_{\pm}=\pm (0.5,0.5)$ and saddle in $(x_1,x_2)=(0,0)$. We
want to consider the non-standard question on how a noise-induced
transition between the two metastable fixed points is achieved in the
most likely way in the small noise limit, $\eps\to0$, if we only
excert stochastic forcing on the $x_1$-component of the
equation. Since the transition necessitates a change in both the $x_1$
and $x_2$ component, but only $x_1$ is fluctuating, we expect a
non-trivial transition trajectory that makes optimal use of the
coupling terms. Additionally, we pick the noise to be multiplicative,
and in particular vanish completely on the line $x_2=x_1+0.2$,
becoming more intense as the solution moves further away from it. The
situation is depicted in Figure~\ref{fig:SDE-2D}~a), where we display
the deterministic dynamics as streamlines, all three fixed points (two
stable and one saddle) as white markers, as well as the noise strength
as shading, with a dashed line covering points where the noise
vanishes. The instanton, as minimiser of the Freidlin-Wentzell
action~\eqref{eq:FW-optim}, for $u_0=u_-$, $u_T=u_+$ and $T=10$, is
depicted as solid line. It is clear that the noise leads the solution
against the flow in a nontrivial and curved way because noise is
available solely on the $x_1$-component, only to use the deterministic
dynamics to approach the saddle. Once crossing the separatrix at the
saddle point, the instanton can relax deterministically into the
opposite stable fixed point. Figure \ref{fig:SDE-2D}~b) shows the
functional Lagrange multiplier $\mu(t)$ as a function of time in
red. In this context, the Lagrange multiplier can be interpreted as the
optimal noise driving the transition. Stochasticity is only needed to
perform the uphill portion of the dynamics, up until the saddle is
reached at approximately $t=7$. After that, the path relaxes
deterministically towards $u_T$ in the proximity to the saddle, and we
have $\mu(t)=0$ for $t>7$. The $x_1(t)$ component is shown alongside
for comparison. We stress that while for a gradient diffusion in
detailed balance, the instanton should follow the flowlines at
every point, the system~\eqref{eq:SDE-2D} breaks detailed balance due
to the multiplicative and degenerate nature of the noise. It this
therefore also not correct that the system's Freidlin-Wentzell
quasipotential is equal to the system's potential.

\subsection{One-dimensional Allen-Cahn model with boundary additive noise}
\label{sec:one-dimens-allen}

\begin{figure}[h!]
  \centering
  \subfloat{
    \includegraphics[width=0.47\linewidth]{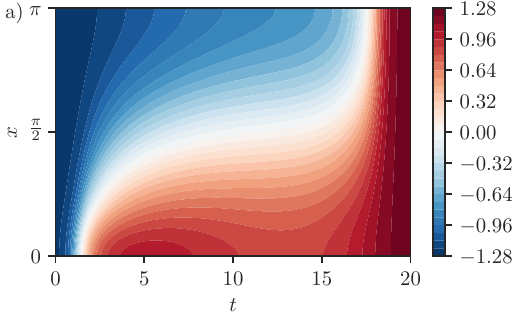}
  }
  \hspace{5mm}
  \subfloat{
        \includegraphics[width=0.47\linewidth]{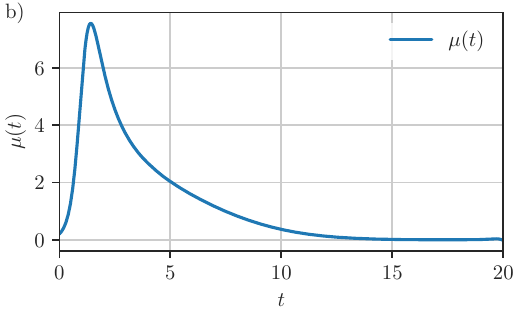}
  }
  \caption{a) Instanton for the transition between
    $u_0\equiv-\sqrt{\alpha}$ to $u_T\equiv\sqrt{\alpha}$ with random Neumann
    boundary condition at $x=0$. First, the noise drags the solution
    towards the saddle by generating an influx at $x=0$. After
    reaching a saddle, at approximately $t=15$, the trajectory
    relaxes deterministically into the other stable solution
    $u_T\equiv\sqrt{\alpha}$. b) The effect of the optimal noise is
    described by the Lagrange multiplier $\mu(t)$, which only acts on
    the portion before reaching the saddle. The results are obtained
    for $\lambda=200$, $\operatorname{F}=\operatorname{I}$ and
    $tol=10^{-3}$.}
  \label{fig:boundary}
\end{figure}

The second example covers the one-dimensional Allen-Cahn model on an
interval $\mathcal{X}=[0,\pi]$. If we would equip the system with
additive Gaussian stochastic forcing, the system would be in detailed
balance with respect to the free energy functional
\begin{equation*}
  E[u] = \int_{\mathcal X}  \left(\tfrac14u^4-\tfrac12\alpha u^2 + (\partial_x u)^2\right)\,\text{d}x\,.
\end{equation*}
Instead, though, we employ stochastic forcing through Neumann
boundary conditions. Concretely, on $x=0$ we consider white Neumann
boundary noise and on $x=\pi$ we set homogeneous Neumann
conditions. Boundary noise has been studied less extensively than
noise across the entire space $\mathcal{X}$. Still,
existence~\cite{da1993evolution},
regularity~\cite{sowers1994multidimensional} and behaviour under small
noise~\cite{freidlin1992reaction} of solutions have been discussed for
different reaction-diffusion equations. For our system,
\begin{align*}
    \left\{\begin{alignedat}{2}
            &\text{d}u^\eps(x,t) = (\partial_x^2 u^\eps (x,t) + \alpha u^\eps (x,t) - u^\eps (x,t)^3) \text{d} t \quad &&, \; x\in\mathcal{X}\;, \; t\geq 0 \;,\\
            &\partial_x u^\eps(0,t) = \sqrt{\eps} \sigma_0 \dot{W}_t &&, \; t\geq 0 \;,\\
            &\partial_x u^\eps(\pi,t) = 0 &&, \; t\geq 0 \;,
    \end{alignedat}\right.
\end{align*}
we follow the analytic results from~\cite{da1993evolution} to obtain
the instanton for the transitions between the two spatially
homogeneous stable fixed points $u_0\equiv-\sqrt{\alpha}$,
$u_T\equiv\sqrt{\alpha}$, $T=20$, for the choice $\alpha=1.5$ and $\sigma_0=0.5$. The
solution of the system solves
\begin{align*}
    u^\eps(t)&= e^{\Delta_\text{N} t} u_0 + \int_0^t e^{\Delta_\text{N} (t-s)} \left(\alpha u^\eps(s)-{u^\eps(s)}^3\right) \text{d}s \\
    &+ \left(\Delta_\text{N}+\operatorname{I}\right) \int_0^t e^{\Delta_\text{N} (t-s)} \operatorname{D} \text{d}W_s \;,
\end{align*}
for any $t\geq0$, $\Delta_\text{N}$ the Laplacian with homogeneous
Neumann boundary conditions in $\mathcal{X}$ and $W$ a Wiener
process. Lastly, $\operatorname{D}:\mathbb{R}\to L^2(\mathcal{X})$
satisfies for any $c_1\in\mathbb{R}$,
\begin{align*}
    \operatorname{D} (c_1)(x)
    = -\frac{\text{cosh}(\pi-x)}{\text{sinh}(\pi)}c_1 \;.
\end{align*}
We also denote as $\operatorname{D}^*:L^2(\mathcal{X})\to\mathbb{R}$
the operator such that, for any $\varphi\in L^2(\mathcal{X})$, it
holds
\begin{align*}
    \operatorname{D}^* (\varphi)
    = -\frac{1}{\text{sinh}(\pi)} \int_0^\pi \text{cosh}(\pi-x) \varphi(x) \text{d}x \;.
\end{align*}
It follows that the covariance operator is
\begin{align*}
    a(x) = a = \left(\Delta_\text{N}+1\right) \operatorname{D} \operatorname{D}^* \left(\Delta_\text{N}+1\right)\;.
\end{align*}
This is all the information needed to compute the large deviation
minimiser for the boundary-noise induced transition between $u_0$ and
$u_T$. Figure \ref{fig:boundary}~a) displays the instanton $\phi$ in space
and time: starting at $u_0\equiv-\sqrt{\alpha}$, the stochastic
forcing generates an influx at $x=0$ that is just enough to push the
system towards the saddle configuration at approximately $t=15$. In
Figure \ref{fig:boundary}~b), we see the Lagrange multiplier $\mu(x,t)$
on $x=0$, where it is analytically not zero. The optimal noise is
concentrated close to $t=0$ and is approximately zero on the
``downhill'' portion for $15<t\le T$, as expected.

\subsection{Spike merging in Gierer-Meinhardt model with degenerate multiplicative noise}
\label{sec:spike-merging-gierer}

\begin{figure}[h!]
  \centering
  \subfloat{
  \includegraphics[width=0.47\linewidth]{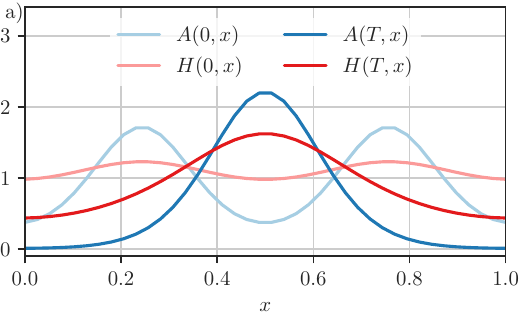}      
  }
  \hspace{5mm}
    \subfloat{
  \includegraphics[width=0.47\linewidth]{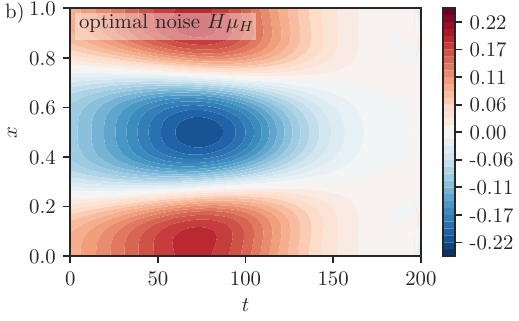}
  }
  \vspace{2mm}
    \subfloat{
  \includegraphics[width=0.47\linewidth]{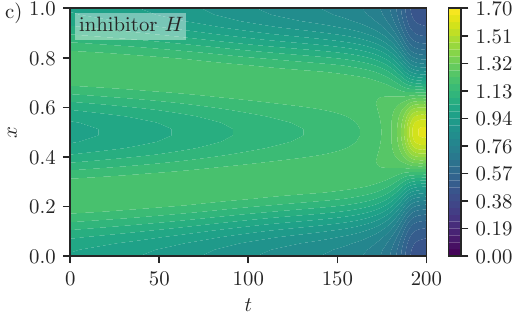}      
  }
  \hspace{5mm}
    \subfloat{
  \includegraphics[width=0.47\linewidth]{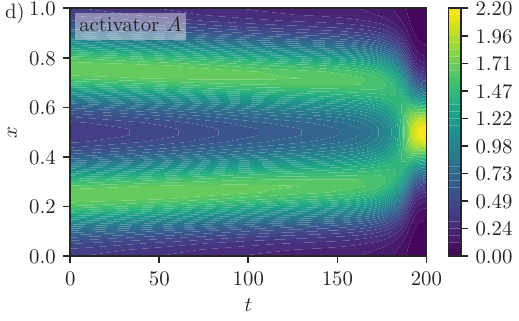}      
  }
  \caption{a) Initial (light) and final (dark) conditions of the
    instanton, transitioning from the stable two-spike configuration
    to the stable one-spike configuration. Here, the inhibitor $H$ is
    depicted in red, and the activator $A$ in blue. b) Optimal noise
    $H \mu_H$ on the inhibitor component. The optimal noise is only
    needed to transiently move the peaks closer together, so that
    afterward, the gap closes deterministically. Figures c) and d) show the
    inhibitor and activator components of the instanton,
    respectively. The two components display the merge of the spikes
    at similar times. The figures have been obtained with the parameters $\lambda=20$, $\operatorname{F}=\operatorname{I}$ and $tol=10^{-4}$.}
    \label{fig:gierer-meinhardt}
\end{figure}

The reaction-diffusion Gierer-Meinhardt
model~\cite{wei2013mathematical} finds applications in biology, such
as on the pattern formation of stripes on seashells. We consider a
one-dimensional version on the interval $\mathcal X=[0,1]$ with homogeneous Neumann boundary
conditions,
\begin{align*}
    \left\{\begin{alignedat}{2}
                &\text{d} u^\eps =
                \begin{pmatrix}
                    \text{d} A \\ \text{d} H
                \end{pmatrix}
                =
                \begin{pmatrix}
                d^2 \partial_x^2 A - A + \frac{A^2}{H}\\
                \frac{1}{\tau} \left( D \partial_x^2 H - H + A^2 \right)
            \end{pmatrix} \text{d} t + \sqrt{\eps} \sigma_0
            \begin{pmatrix}
                0 \\ H
            \end{pmatrix} \text{d} W_t &&,\\
            &\partial_x A(0,t) = \partial_x A(1,t) = 0 &&, \; t\geq 0 \;,\\
            &\partial_x H(0,t) = \partial_x H(1,t) = 0 &&, \; t\geq 0 \;.
    \end{alignedat}\right.
\end{align*}
The Gierer-Meinhardt model is known to display steady solutions
characterised by spikes on the variables $A\geq0$, the so-called
``activator'', and $H\geq0$, the ``inhibitor''. The number of spikes
present in a stable solution depends on the diffusivity constants
$d,D>0$. Moreover, the stability of the flat solution $A=H\equiv1$,
depends on the timescale
$\tau>0$. Following~\cite{antwi2019global,kelkel2010stochastic,winter2016dynamics},
we consider the effect of degenerate multiplicative noise that forces
the component $H$ only, leaving $A$ to act solely under the influence
of $H$. In the figures to follow, we discuss the instanton that
describes the merge of two spikes in time $T=200$ and for $d=0.06$,
$D=0.04$, $\tau=0.5$ and $\sigma_0=0.5$. In Figure \ref{fig:gierer-meinhardt}~a), we
display the initial conditions $u_0$, light, and the final conditions
$u_T$, dark. The steady solutions $u_0$ and $u_T$ are obtained
numerically following~\cite[Chapter 2.1]{wei2013mathematical} and
their stability is proven by~\cite[Remark
  4.4]{wei2013mathematical}. We are interested in the effect of
noise-induced \textit{spike merging}, i.e.~the stochasticity
transforming a pair of spikes into a single spike.

In Figure \ref{fig:gierer-meinhardt}~a), the forced component $H$ is
shown in red, and the activator $A$ is indicated in blue, for the
initial conditions (light color) and final conditions (dark color). In
Figure \ref{fig:gierer-meinhardt}~b), the optimal noise is displayed
as $H \mu_H$, for $\mu_H$ the Lagrange multiplier associated with the
differential equation of $H$. It is apparent that the noise
prioritises shifting the spikes closer together and, afterward, when a
critical distance is achieved at roughly $t=150$, the instanton
approaches $u_T$ deterministically. In Figure
\ref{fig:gierer-meinhardt}~c) and in Figure
\ref{fig:gierer-meinhardt}~d), the components $H$ and $A$ are shown,
respectively. As the instanton approaches $u_T$ in a deterministic
manner, the merge of the two spikes occurs
simultaneously in the inhibitor and the activator components.

\subsection{Pulse initiation in FitzHugh-Nagumo model with degenerate additive noise}
\label{sec:pulse-init-fitzh}

\begin{figure}[h!]
  \centering
  \subfloat{
  \includegraphics[width=0.47\linewidth]{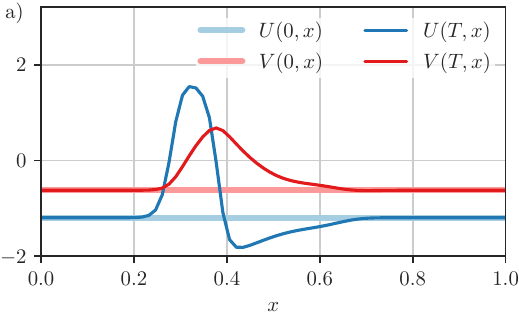}      
  }
  \hspace{5mm}
  \subfloat{
  \includegraphics[width=0.47\linewidth]{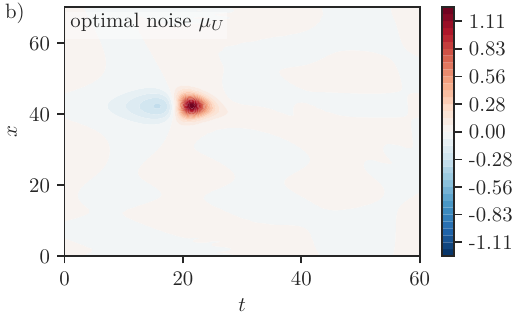}    
  }
  \vspace{2mm}
  \subfloat{
  \includegraphics[width=0.47\linewidth]{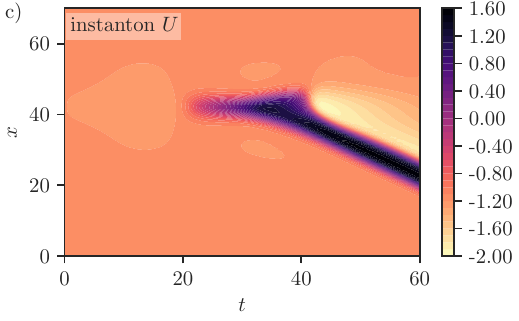}     
  }
  \hspace{5mm}
  \subfloat{
  \includegraphics[width=0.47\linewidth]{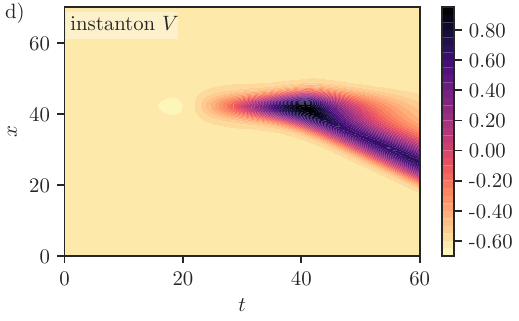}     
  }
  \caption{a) The initial, $u_0$, and final, $u_T$, conditions
    in light and dark lines, respectively. The forced component, $U$,
    is displayed in blue, and $V$ is indicated in red. b) The Lagrange
    multiplier is shown and indicates the effect of the noise in
    the pulse initiation. Figures c) and d) show the pulse initiation on the
    $U$ and $V$ components, respectively. The figures are obtained for
    parameters $\lambda=0.5$ and $tol = 5\;10^{-4}$ and with periodic
    boundary conditions. The filter on the final condition is
    $\operatorname{F}=\operatorname{I}$.}
  \label{fig:FHN}
\end{figure}

The reaction-diffusion FitzHugh-Nagumo model is often employed in the
simulation of electric impulses through nerve
axons~\cite{ermentrout2010mathematical,izhikevich2007dynamical}. It
presents a behaviour qualitatively similar to the Hodgkin-Huxley
model~\cite{hodgkin1952quantitative}, despite being composed of
significantly simpler equations. The reaction-diffusion model on the
real line is characterised by traveling pulse solutions, whose
properties have been extensively
studied~\cite{guckenheimer2012homoclinic,sandstede1998stability}. In
the current subsection, we construct an instanton for the model
\begin{align}
    \text{d} u^\eps =
\begin{pmatrix}
        \text{d}U \\ \text{d}V
    \end{pmatrix}
    =
    \begin{pmatrix}
        \nu_1 \partial_{x}^2 U + U - U^3 - V \\
        \nu_2 \partial_{x}^2 V + \delta \left( U - \gamma_1 V + \gamma_2 \right)
    \end{pmatrix} \text{d} t
     + \sqrt{\eps} \sigma_0
            \begin{pmatrix}
                1 \\ 0
            \end{pmatrix} \text{d} W_t \;,
            \label{eq:FHN}
\end{align}
the reaction-diffusion FitzHugh-Nagumo
model with additive noise on variable $U$, following the example
of~\cite{eichinger2022multiscale}. We associate the component $U$ to
the electric potential and $V$ with a recovery variable. The instanton
describes the noise-induced initiation of a pulse, an event
well-studied~\cite{idris2008initiation}. Note that in~\eqref{eq:FHN},
only the electric potential is subject to (additive Gaussian)
stochastic noise, while the recovery variable is left unforced. The
parameters are set as $T=60$, $\nu_1=1$, $\nu_2=0.1$, $\delta=0.08$,
$\gamma_1=0.8$, $\gamma_2=0.7$ and $\sigma_0=0.5$. The spatially homogeneous initial
condition $u_0\approx (-1.19941,-0.62426)^T$ (absence of a pulse) and
final condition $u_T$ (pulse present) are chosen from~\cite[Example
  2.4]{beyn2018computation}. They are displayed in Figure
\ref{fig:FHN}~a) in light and dark lines, respectively. The forced
component $U$ is indicated in blue, and the term $V$ is shown in
red. The stability of $u_0$ can be easily computed, and the stability
of the traveling wave of frame $u_T$ is obtained in~\cite[Example
  3.7]{beyn2018computation}. In Figure \ref{fig:FHN}~b), the Lagrange
multiplier, $\mu_U$, associated with the first equation in the model,
is displayed. The noise provides first a small negative push to the
flat solution $u_0$ in a concentrated region and, secondly, initiates
the pulse with a larger input. Then, it directs the pulse in the same
direction along which $u_T$ travels. In fact, the traveling wave
mirrored in space to the pulse associated to $u_T$ is also a stable
solution of the system. Lastly, the instanton converges
deterministically to $u_T$. In Figure \ref{fig:FHN}~c) and in Figure
\ref{fig:FHN}~d), the perturbed component $U$ and the component $V$
are shown in contour plots, respectively. The creation of the pulse
appears to be close to simultaneous in the terms. Furthermore, the
negative bell in the tail of the pulse in the term $U$ arises in a
deterministic manner. Note that the pulse solution is only a fixed
point in a reference frame of its movement speed, while our equations
are defined in the laboratory reference frame. For this reason, the
instanton, as solution to the optimisation problem, automatically
initiates the pulse at a sufficient distance to allow it to travel to
its pre-chosen endpoint.

\subsection{Puff splitting in Barkley model with degenerate multiplicative noise}
\label{sec:puff-splitt-barkl}

\begin{figure}[h!]
  \centering
  \subfloat{
  \includegraphics[width=0.47\linewidth]{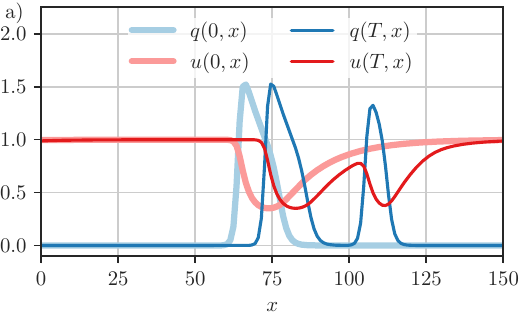}      
  }
  \hspace{5mm}
  \subfloat{
  \includegraphics[width=0.47\linewidth]{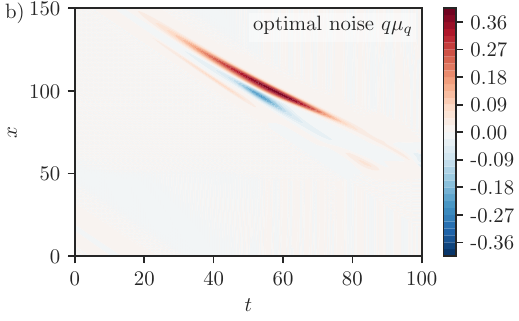}      
  }
  \vspace{2mm}
  \subfloat{
  \includegraphics[width=0.47\linewidth]{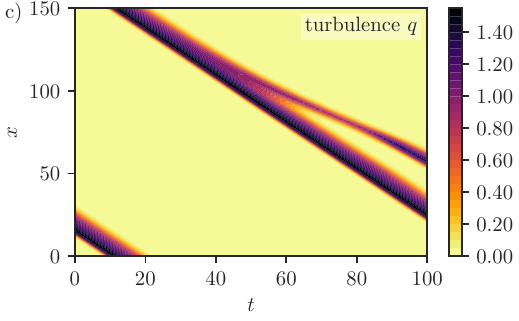} 
  }
  \hspace{5mm}
  \subfloat{
  \includegraphics[width=0.47\linewidth]{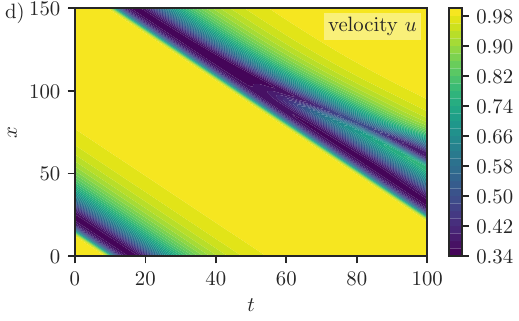} 
  }
  \caption{a) The initial, light, and
    final, dark, conditions of the instanton. The enforced final
    condition, on $\mathcal{X}_1=[50,70]$, is shown in green. The blue
    lines are associated with the forced component $q$ and the red
    lines to $u$. The state $u_0$ is a stable solution of the model
    but $u_T$ is not the frame of a traveling wave. Figure b) displays $q \mu_q$, which indicates the optimal noise that defines the instanton. Figures c) and d) show the puff splitting event for the components $q$ and $u$, respectively. The figures are obtained
    for parameters $\lambda = 200$ and $tol = 10^{-2}$.}
    \label{fig:Barkley}
\end{figure}

As final, and from a numerical perspective most complex, last example,
we take the Barkley model for the evolution of turbulent puffs in pipe
flow. It is defined by coupled SPDEs on the real line with degenerate
multiplicative noise,
\begin{comment}
\tiny
\begin{align*}
    \dot{u}^\eps =
    \begin{pmatrix}
        \dot{q} \\ \dot{u}
    \end{pmatrix}
    =
    \begin{pmatrix}
        D \partial_{x}^2 q + \left(\xi-u\right) \partial_x q + q \left( r + u - U_0 - (r+\delta) (q-1)^2 \right) \\
        - u \partial_x u + \eps_1 ( U_0 - u ) + \eps_2 q (\bar{U} - u )
    \end{pmatrix}
     + \sqrt{\eps} \sigma
            \begin{pmatrix}
                q \\ 0
            \end{pmatrix} \eta
\end{align*}
\normalsize
\end{comment}
\begin{align*}
    \text{d} u^\eps \!=\!
    \begin{pmatrix}
        \text{d} q \\ \text{d} u
    \end{pmatrix}
    \!=\!
    \begin{pmatrix}
        D \partial_{x}^2 q + \xi_2  \partial_x q + q \left( u + r - 1 - (r+\delta) (q-1)^2 \right) \\
        \left(-\xi_1 + \xi_2 \right) \partial_x u + \eps_1 ( 1 - u ) - \eps_2 u q
    \end{pmatrix} \text{d} t
     + \sqrt{\eps} \sigma_0
            \begin{pmatrix}
                q \\ 0
            \end{pmatrix} \text{d} W_t \;.
\end{align*}
The system describes the evolution of turbulence in shear flows
through a long pipe~\cite{barkley2011modeling}, extended in
$x$-direction. The variable $q\geq0$ indicates the turbulent kinetic
energy, as difference of the transverse velocity components to the
laminar background profile. The variable $u$ in turn describes the
centerline velocity. The model is normalised such that the centerline
velocity is $0$ in the presence of strong turbulence, and $1$ for the
laminar flow. The variable $r>0$ represents the fluid Reynolds
number. Phenomenologically, this model always allows for the spatially
homogeneous laminar flow $(q,u)=(0,1)$ to be stable at any Reynolds
number, representing the fact that in transitional pipe flow the
laminar solution of the Navier-Stokes equation remains linearly
stable. The stochastic forcing models the chaotic nature of turbulent
flow, and is hence chosen to be acting on $q$ only, and proportional
to $q$ itself. Consequently, the laminar solution exhibits no noise at
all and is hence an absorbing state. At sufficient Reynolds numbers,
the system additionally exhibits turbulent puff solutions in the form
of localised travelling packets with $q(x)>0$, which similarly occur
in actual pipe flow. They are long-lived in nature, and linearly
stable in the model. Crucially, stochasticity may force the turbulent
puff to either decay into laminar flow, or alternatively split into
two independent, separated turbulent puffs. Above a critical Reynolds
number $r_c$, puff splitting dominates puff decay, and turbulence
proliferates in the pipe~\cite{barkley2011modeling,
  barkley2016theoretical, gome2022extreme}.

Our aim is to capture the puff splitting
process~\cite{frishman-grafke:2022-a} by computing the split
instanton, the noise-induced transition from a single puff into
two. Here, the noise generates a second puff splitting off from an
existing one, on the interval $[0,150]$ with periodic boundary
conditions. Note that, in contrast to the model in subsection~\ref{sec:pulse-init-fitzh}, a puff cannot be created in a
completely laminar region, since there the stochastic forcing is
necessarily identically zero. This makes this model particularly
difficult to handle, as the noise is not only inactive on the whole
$u$-field, but dynamically inactive for $q$ in most of the domain as
well.

The parameters are taken as $D=0.5$ for the diffusion, $r=0.6$,
$\delta = 0.1$ and $\sigma_0=0.5$. In order to describe the centerline velocity, we have
chosen the values $\eps_1=0.1$ and $\eps_2=0.2$. Lastly, the
parameters $\xi_1=\xi_2=0.8$ refer to the advection speed and the
moving reference frame. The starting condition $u_0$ is a puff,
obtained as a converging solution from an initial bell-shaped
state. The initial condition is shown in Figure \ref{fig:Barkley} a)
in a light blue line for the perturbed component $q$ and a light red
line for $u$. The final condition is obtained as follows. The Barkley
model is considered with the initial condition $u_0$ at $t=0$. We
label as $u_1=(q_1,u_1)^T$ the state at $t=70$ and as $u_2$ the state
at $t=100$. We define as $\mathfrak{R}_c$ the rotation operator in the
right direction on the interval $[0,150]$ with periodic boundary
conditions for the value $c$. To construct $u_T$, we observe the Barkley model with
initial conditions in
$u_3=(q_3,u_3)^T=(q_1+\mathfrak{R}_{\frac{150}{7}} q_1,
u_1+\mathfrak{R}_{\frac{150}{7}} u_1-1)^T$ at $t=0$. The final
condition $u_T$ is defined as the state of the solution of the system
at time $t=30$. Its shape corresponds to a leading and a second
well-separated trailing puff, as shown in Figure \ref{fig:Barkley} a)
in dark lines. The component $q$ is shown in blue, and $u$ is
indicated in red. The state $u_3$ is chosen to obtain a realistic
sustainable double puff state in $u_T$. Furthermore, at time $T=100$,
the position of the higher puff in $u_T$ is qualitatively similar to
the location of the puff shown in the state $u_2$ on the interval with
periodic boundary conditions. A filter is chosen in order to capture
the shape of the second puff in $u_T$. Therefore
$\operatorname{F}=\mathbbm{1}_{\mathcal{X}_1}$, the multiplication
operator for the indicator function with support on
$\mathcal{X}_1=[50,70]$. The solid green lines indicate the enforced
shape in the penalty method.

In Figure \ref{fig:Barkley} b), the optimal noise is described by $q
\mu_q$, for $\mu_q$ that indicates the Lagrange multiplier for the
first differential equation in the system. The conjugate variable can
be interpreted to cause the following: first, it elongates the tail of
the original puff into a wider object; second, for $t\in[40,60]$, it
generates a laminar gap to cut the elongated puff into two, deep enough
to result in a split; afterward, weak final stochastic adjustments
are made to meet the final locations of the two puffs. The
adjustments are further justified by the fact that $u_T$ is not the
frame of a steady traveling wave. The instanton describing a puff
splitting event is shown in Figure \ref{fig:Barkley} c) and in Figure
\ref{fig:Barkley} d), which display the components $q$ and $u$,
respectively, in contour plots. Similarly to the previous examples,
the splitting event appears to happen simultaneously on both terms,
despite the degeneracy of the noise.

\section{Conclusion}

We introduce a method to compute instantons, large deviation
minimisers, for metastable stochastic partial differential equations
in the presence of degenerate noise that are inaccessible to
computation by existing methods. We demonstrate the applicability of
the method for various example problems of increasing complexity,
culminating in the computation of the turbulent puff splitting
instanton in the Barkley model, a system that describes the behaviour
of regions of turbulence in a long pipe.

Further applications to the algorithm are spike merging in the
Gierer-Meinhardt model, with multiplicative noise on the inhibitor
component, and pulse creation in the FitzHugh-Nagumo model, with
additive noise on the electric potential component. The implementation
of the algorithm enables the construction of the most likely paths
under degenerate noise and boundary noise, for different types of
models.

\section*{Acknowledgments}

The authors express their gratitude to Christian Kuehn and Freddy
Bouchet for the many valuable discussions, especially regarding the
Barkley model.

\newpage
\printbibliography

@book{freidlin1998random,
  title={Random Perturbations of Dynamical Systems},
  author={Freidlin, Mark Iosifovich and Wentzell, Alexander D},
  year={1998},
  publisher={Springer}
}

@book{da2014stochastic,
  title={Stochastic equations in infinite dimensions},
  author={Da Prato, Giuseppe and Zabczyk, Jerzy},
  volume={152},
  year={2014},
  publisher={Cambridge university press}
}

@article{rolland2016computing,
  title={Computing transition rates for the 1-D stochastic Ginzburg--Landau--Allen--Cahn equation for finite-amplitude noise with a rare event algorithm},
  author={Rolland, Joran and Bouchet, Freddy and Simonnet, Eric},
  journal={J. Stat. Phys.},
  volume={162},
  pages={277--311},
  year={2016},
  publisher={Springer}
}

@article{grafke2019numerical,
  title={Numerical computation of rare events via large deviation theory},
  author={Grafke, Tobias and Vanden-Eijnden, Eric},
  journal={Chaos},
  volume={29},
  number={6},
  year={2019},
  publisher={AIP Publishing}
}

@article{berglund2023stochastic,
  title={Stochastic resonance in stochastic PDEs},
  author={Berglund, Nils and Nader, Rita},
  journal={Stoch. Partial Differ. Equ.},
  volume={11},
  number={1},
  pages={348--387},
  year={2023},
  publisher={Springer}
}

@article{bernuzzi2023bifurcations,
  title={Bifurcations and early-warning signs for spdes with spatial heterogeneity},
  author={Bernuzzi, Paolo and Kuehn, Christian},
  journal={J. Dyn. Differ. Equ.},
  pages={1--45},
  year={2023},
  publisher={Springer}
}

@book{nocedal1999numerical,
  title={Numerical optimization},
  author={Nocedal, Jorge and Wright, Stephen J},
  year={1999},
  publisher={Springer}
}

@article{da1993evolution,
  title={Evolution equations with white-noise boundary conditions},
  author={Da Prato, Giuseppe},
  journal={Stochastics},
  volume={42},
  number={3-4},
  pages={167--182},
  year={1993},
  publisher={Taylor \& Francis}
}

@article{freidlin1992reaction,
  title={Reaction-diffusion equations with randomly perturbed boundary conditions},
  author={Freidlin, Mark I and Wentzell, Alexander D},
  journal={Ann. Probab.},
  pages={963--986},
  year={1992},
  publisher={JSTOR}
}

@article{sowers1994multidimensional,
  title={Multidimensional reaction-diffusion equations with white noise boundary perturbations},
  author={Sowers, Richard B},
  journal={Ann. Probab.},
  pages={2071--2121},
  year={1994},
  publisher={JSTOR}
}

@book{wei2013mathematical,
  title={Mathematical aspects of pattern formation in biological systems},
  author={Wei, Juncheng and Winter, Matthias},
  volume={189},
  year={2013},
  publisher={Springer Science \& Business Media}
}

@article{kelkel2010stochastic,
  title={On a stochastic reaction--diffusion system modeling pattern formation on seashells},
  author={Kelkel, Jan and Surulescu, Christina},
  journal={J. Math. Biol.},
  volume={60},
  pages={765--796},
  year={2010},
  publisher={Springer}
}

@article{winter2016dynamics,
  title={The dynamics of the stochastic shadow Gierer--Meinhardt system},
  author={Winter, Matthias and Xu, Lihu and Zhai, Jianliang and Zhang, Tusheng},
  journal={J. Differ. Equ.},
  volume={260},
  number={1},
  pages={84--114},
  year={2016},
  publisher={Elsevier}
}

@article{antwi2019global,
  title={Global analysis of the shadow Gierer-Meinhardt system with general linear boundary conditions in a random environment},
  author={Antwi-Fordjour, Kwadwo and Kim, Seonguk and Nkashama, Marius},
  journal={arXiv preprint arXiv:1907.10122},
  year={2019}
}

@book{ermentrout2010mathematical,
  title={Mathematical foundations of neuroscience},
  author={Ermentrout, Bard and Terman, David Hillel},
  volume={35},
  year={2010},
  publisher={Springer}
}

@book{izhikevich2007dynamical,
  title={Dynamical systems in neuroscience},
  author={Izhikevich, Eugene M},
  year={2007},
  publisher={MIT press}
}

@article{hodgkin1952quantitative,
  title={A quantitative description of membrane current and its application to conduction and excitation in nerve},
  author={Hodgkin, Alan L and Huxley, Andrew F},
  journal={J. Physiol.},
  volume={117},
  number={4},
  pages={500},
  year={1952},
  publisher={Wiley}
}

@article{guckenheimer2012homoclinic,
  title={Homoclinic orbits of the FitzHugh-Nagumo equation: The singular-limit},
  author={Guckenheimer, John and Kuehn, Christian},
  journal={arXiv preprint arXiv:1201.5901},
  year={2012}
}

@article{sandstede1998stability,
  title={Stability of multiple-pulse solutions},
  author={Sandstede, Bj{\"o}rn},
  journal={Trans. Am. Math. Soc.},
  volume={350},
  number={2},
  pages={429--472},
  year={1998}
}

@article{eichinger2022multiscale,
  title={Multiscale analysis for traveling-pulse solutions to the stochastic FitzHugh--Nagumo equations},
  author={Eichinger, Katharina and Gnann, Manuel V and Kuehn, Christian},
  journal={Ann. Appl. Probab.},
  volume={32},
  number={5},
  pages={3229--3282},
  year={2022},
  publisher={Institute of Mathematical Statistics}
}

@book{idris2008initiation,
  title={Initiation of excitation waves},
  author={Idris, Ibrahim},
  year={2008},
  publisher={The University of Liverpool (United Kingdom)}
}

@article{beyn2018computation,
  title={Computation and stability of traveling waves in second order evolution equations},
  author={Beyn, W-J and Otten, Denny and Rottmann-Matthes, Jens},
  journal={SIAM J. Numer. Anal.},
  volume={56},
  number={3},
  pages={1786--1817},
  year={2018},
  publisher={SIAM}
}

@inproceedings{barkley2011modeling,
  title={Modeling the transition to turbulence in shear flows},
  author={Barkley, Dwight},
  booktitle={J. Phys. Conf. Ser.},
  volume={318},
  number={3},
  pages={032001},
  year={2011},
  organization={IOP Publishing}
}

@article{barkley2016theoretical,
  title={Theoretical perspective on the route to turbulence in a pipe},
  author={Barkley, Dwight},
  journal={J. Fluid Mech.},
  volume={803},
  pages={P1},
  year={2016},
  publisher={Cambridge University Press}
}

@article{plessix2006review,
  title={A review of the adjoint-state method for computing the gradient of a functional with geophysical applications},
  author={Plessix, R-E},
  journal={Geophys. J. Int.},
  volume={167},
  number={2},
  pages={495--503},
  year={2006},
  publisher={Blackwell Publishing Ltd Oxford, UK}
}

@article{gome2022extreme,
  title={Extreme events in transitional turbulence},
  author={Gom{\'e}, S{\'e}bastien and Tuckerman, Laurette S and Barkley, Dwight},
  journal={Philos. Trans. R. Soc. A},
  volume={380},
  number={2226},
  pages={20210036},
  year={2022},
  publisher={The Royal Society}
}

@article{borovykh2023privacy,
  title={Privacy Risk for anisotropic Langevin dynamics using relative entropy bounds},
  author={Borovykh, Anastasia and Kantas, Nikolas and Parpas, Panos and Pavliotis, Greg},
  journal={arXiv preprint arXiv:2302.00766},
  year={2023}
}

@article{blomker2005modulation,
  title={Modulation equations: Stochastic bifurcation in large domains},
  author={Bl{\"o}mker, Dirk and Hairer, M and Pavliotis, GA},
  journal={Commun. Math. Phys.},
  volume={258},
  pages={479--512},
  year={2005},
  publisher={Springer}
}

@article{khoo2019solving,
  title={Solving for high-dimensional committor functions using artificial neural networks},
  author={Khoo, Yuehaw and Lu, Jianfeng and Ying, Lexing},
  journal={Res. Math. Sci.},
  volume={6},
  pages={1--13},
  year={2019},
  publisher={Springer}
}

@article{liu2015efficient,
  title={Efficient rare event simulation for failure problems in random media},
  author={Liu, Jingchen and Lu, Jianfeng and Zhou, Xiang},
  journal={SIAM J. Sci. Comput.},
  volume={37},
  number={2},
  pages={A609--A624},
  year={2015},
  publisher={SIAM}
}

@article{bisewski2019rare,
  title={Rare event simulation for steady-state probabilities via recurrency cycles},
  author={Bisewski, Krzysztof and Crommelin, Daan and Mandjes, Michel},
  journal={Chaos},
  volume={29},
  number={3},
  year={2019},
  publisher={AIP Publishing}
}

@article{crommelin2015stochastic,
  title={Stochastic and statistical methods in climate, atmosphere, and ocean science},
  author={Crommelin, Daan and Khouider, Boualem and others},
  year={2015},
  publisher={HeidelbergSpringer Reference}
}

@article{xu2014rare,
  title={Rare-Event Simulation for the Stochastic Korteweg--de Vries Equation},
  author={Xu, Gongjun and Lin, Guang and Liu, Jingchen},
  journal={SIAM/ASA J. Uncertain. Quantificat.},
  volume={2},
  number={1},
  pages={698--716},
  year={2014},
  publisher={SIAM}
}

@article{ragone2021rare,
  title={Rare event algorithm study of extreme warm summers and heatwaves over Europe},
  author={Ragone, Francesco and Bouchet, Freddy},
  journal={Geophys. Res. Lett.},
  volume={48},
  number={12},
  pages={e2020GL091197},
  year={2021},
  publisher={Wiley Online Library}
}

@article{ragone2018computation,
  title={Computation of extreme heat waves in climate models using a large deviation algorithm},
  author={Ragone, Francesco and Wouters, Jeroen and Bouchet, Freddy},
  journal={Proc. Natl. Acad. Sci. U.S.A.},
  volume={115},
  number={1},
  pages={24--29},
  year={2018},
  publisher={National Acad Sciences}
}

@article{breden2019rigorous,
  title={Rigorous validation of stochastic transition paths},
  author={Breden, Maxime and Kuehn, Christian},
  journal={J. Math. Pures Appl.},
  volume={131},
  pages={88--129},
  year={2019},
  publisher={Elsevier}
}

@article{blomker2012numerical,
  title={Numerical study of amplitude equations for SPDEs with degenerate forcing},
  author={Bl{\"o}mker, Dirk and Mohammed, Wael W and Nolde, Christian and W{\"o}hrl, Franz},
  journal={Int. J. Comput. Math.},
  volume={89},
  number={18},
  pages={2499--2516},
  year={2012},
  publisher={Taylor \& Francis}
}

@article{dong2018ergodicity,
  title={Ergodicity of the 2D Navier-Stokes equations with degenerate multiplicative noise},
  author={Dong, Zhao and Peng, Xu-hui},
  journal={Acta Math. Appl. Sin.},
  volume={34},
  number={1},
  pages={97--118},
  year={2018},
  publisher={Springer}
}

@article{baars2017continuation,
  title={Continuation of probability density functions using a generalized Lyapunov approach},
  author={Baars, Sven and Viebahn, JP and Mulder, Thomas E and Kuehn, Christian and Wubs, Fred W and Dijkstra, Henk A},
  journal={J. Comput. Phys.},
  volume={336},
  pages={627--643},
  year={2017},
  publisher={Elsevier}
}

@article{bernuzzi2024warning,
  title={Warning Signs for Boundary Noise and their Application to an Ocean Boussinesq Model},
  author={Bernuzzi, Paolo and Dijkstra, Henk A and Kuehn, Christian},
  journal={arXiv preprint arXiv:2405.13550},
  year={2024}
}

@book{rubino2009rare,
  title={Rare event simulation using Monte Carlo methods},
  author={Rubino, Gerardo and Tuffin, Bruno and others},
  volume={73},
  year={2009},
  publisher={Wiley Online Library}
}

@Article{	  alqahtani-grafke:2021,
  title		= {Instantons for rare events in heavy-tailed distributions},
  volume	= {54},
  issn		= {1751-8121},
  doi		= {10.1088/1751-8121/abe67b},
  number	= {17},
  journal	= {J. Phys. A},
  author	= {Alqahtani, Mnerh and Grafke, Tobias},
  month		= apr,
  year		= {2021},
  pages		= {175001}
}

@Article{	  costeniuc-ellis-touchette-etal:2005,
  title		= {The {Generalized} {Canonical} {Ensemble} and {Its}
		  {Universal} {Equivalence} with the {Microcanonical}
		  {Ensemble}},
  volume	= {119},
  issn		= {1572-9613},
  doi		= {10.1007/s10955-005-4407-0},
  number	= {5},
  journal	= {J. Stat. Phys.},
  author	= {Costeniuc, Marius and Ellis, Richard S. and Touchette,
		  Hugo and Turkington, Bruce},
  month		= jun,
  year		= {2005},
  pages		= {1283--1329}
}

@Article{	  e-ren-vanden-eijnden:2004,
  title		= {Minimum action method for the study of rare events},
  volume	= {57},
  issn		= {1097-0312},
  doi		= {10.1002/cpa.20005},
  number	= {5},
  journal	= {Commun. Pure Appl. Math.},
  author	= {E, Weinan and Ren, Weiqing and Vanden-Eijnden, Eric},
  month		= may,
  year		= {2004},
  pages		= {637--656}
}

@Article{	  grafke-grauer-schaefer-etal:2014,
  title		= {Arclength {Parametrized} {Hamilton}'s {Equations} for the
		  {Calculation} of {Instantons}},
  volume	= {12},
  issn		= {1540-3459},
  doi		= {10.1137/130939158},
  number	= {2},
  journal	= {Multiscale Model. Simul.},
  author	= {Grafke, T. and Grauer, R. and Sch\"afer, T. and
		  Vanden-Eijnden, E.},
  month		= jan,
  year		= {2014},
  pages		= {566--580}
}

@InCollection{	  grafke-schaefer-vanden-eijnden:2017,
  series	= {Fields {Institute} {Communications}},
  title		= {Long {Term} {Effects} of {Small} {Random} {Perturbations}
		  on {Dynamical} {Systems}: {Theoretical} and {Computational}
		  {Tools}},
  isbn		= {978-1-4939-6968-5},
%  isbn		= {978-1-4939-6968-5 978-1-4939-6969-2},
  shorttitle	= {Long {Term} {Effects} of {Small} {Random} {Perturbations}
		  on {Dynamical} {Systems}},
  booktitle	= {Recent {Progress} and {Modern} {Challenges} in {Applied}
		  {Mathematics}, {Modeling} and {Computational} {Science}},
  publisher	= {Springer, New York, NY},
  author	= {Grafke, Tobias and Sch\"afer, Tobias and Vanden-Eijnden,
		  Eric},
  year		= {2017},
  doi		= {10.1007/978-1-4939-6969-2\_2},
  pages		= {17--55},
  cmt_editor	= {Melnik, Roderick and Makarov, Roman and Belair, Jacques}
}

@Article{	  heymann-vanden-eijnden:2008-a,
  author	= {Heymann, M. and Vanden-Eijnden, E.},
  journal	= {Phys. Rev. Lett.},
  number	= {14},
  pages		= {140601},
  publisher	= {APS},
  title		= {Pathways of maximum likelihood for rare events in
		  nonequilibrium systems: application to nucleation in the
		  presence of shear},
  volume	= {100},
  year		= {2008}
}

@Article{	  vanden-eijnden-heymann:2008,
  author	= {E. Vanden-Eijnden and M. Heymann},
  title		= {The geometric minimum action method for computing minimum
		  energy paths},
  journal	= {Jour. Chem. Phys.},
  year		= {2008},
  volume	= {128},
  pages		= {061103}
}

@Article{	  zhou-ren-e:2008,
  title		= {Adaptive minimum action method for the study of rare
		  events},
  author	= {Zhou, Xiang and Ren, Weiqing and E, Weinan},
  journal	= {J. Chem. Phys.},
  volume	= {128},
  pages		= {104111},
  year		= {2008}
}

@Article{	  frishman-grafke:2022-a,
  title		= {Mechanism for turbulence proliferation in subcritical
		  flows},
  volume	= {478},
  doi		= {10.1098/rspa.2022.0218},
  number	= {2265},
  journal	= {Proc. R. Soc. A},
  author	= {Frishman, Anna and Grafke, Tobias},
  month		= sep,
  year		= {2022},
  pages		= {20220218}
}

@article{lenton2013origin,
  title={On the origin of planetary-scale tipping points},
  author={Lenton, Timothy M and Williams, Hywel TP},
  journal={Trends Ecol. Evol.},
  volume={28},
  number={7},
  pages={380--382},
  year={2013},
  publisher={Elsevier}
}

@article{loriani2023tipping,
  title={Tipping points in ocean and atmosphere circulations},
  author={Loriani, Sina and Aksenov, Yevgeny and Armstrong McKay, David and Bala, Govindasamy and Born, Andreas and Chiessi, Cristiano M and Dijkstra, Henk and Donges, Jonathan F and Drijfhout, Sybren and England, Matthew H and others},
  journal={EGUsphere},
  volume={2023},
  pages={1--62},
  year={2023},
  publisher={Copernicus Publications G{\"o}ttingen, Germany}
}

@article{dakos2024tipping,
  title={Tipping point detection and early warnings in climate, ecological, and human systems},
  author={Dakos, Vasilis and Boulton, Chris A and Buxton, Joshua E and Abrams, Jesse F and Arellano-Nava, Beatriz and Armstrong McKay, David I and Bathiany, Sebastian and Blaschke, Lana and Boers, Niklas and Dylewsky, Daniel and others},
  journal={Earth Syst. Dyn.},
  volume={15},
  number={4},
  pages={1117--1135},
  year={2024},
  publisher={Copernicus Publications G{\"o}ttingen, Germany}
}

@Article{	  cerou-guyader-rousset:2019,
  title		= {Adaptive multilevel splitting: {Historical} perspective
		  and recent results},
  volume	= {29},
  issn		= {1054-1500},
  shorttitle	= {Adaptive multilevel splitting},
  doi		= {10.1063/1.5082247},
  number	= {4},
  journal	= {Chaos},
  author	= {C\'erou, Fr\'ed\'eric and Guyader, Arnaud and Rousset,
		  Mathias},
  month		= apr,
  year		= {2019},
  pages		= {043108}
}

@Article{	  cerou-guyader:2007,
  title		= {Adaptive multilevel splitting for rare event analysis},
  volume	= {25},
  issn		= {0736-2994},
  doi		= {10.1080/07362990601139628},
  number	= {2},
  journal	= {Stochastic Anal. Appl.},
  author	= {Cerou, Frederic and Guyader, Arnaud},
  year		= {2007},
  pages		= {417--443}
}

@Article{	  lestang-ragone-brehier-etal:2018,
  title		= {Computing return times or return periods with rare event
		  algorithms},
  volume	= {2018},
  issn		= {1742-5468},
  doi		= {10.1088/1742-5468/aab856},
  number	= {4},
  journal	= {J. Stat. Mech: Theory Exp.},
  author	= {Lestang, Thibault and Ragone, Francesco and Br\'ehier,
		  Charles-Edouard and Herbert, Corentin and Bouchet, Freddy},
  month		= apr,
  year		= {2018},
  pages		= {043213}
}

@Article{	  simonnet:2023,
  title		= {Computing non-equilibrium trajectories by a deep learning
		  approach},
  volume	= {491},
  issn		= {0021-9991},
  doi		= {10.1016/j.jcp.2023.112349},
  journal	= {J. Comput. Phys.},
  author	= {Simonnet, Eric},
  month		= oct,
  year		= {2023},
  pages		= {112349}
}

@Article{	  zakine-vanden-eijnden:2023,
  title		= {Minimum-{Action} {Method} for {Nonequilibrium} {Phase}
		  {Transitions}},
  volume	= {13},
  issn		= {2160-3308},
  doi		= {10.1103/PhysRevX.13.041044},
  number	= {4},
  journal	= {Phys. Rev. X},
  author	= {Zakine, Ruben and Vanden-Eijnden, Eric},
  month		= dec,
  year		= {2023},
  pages		= {041044}
}

\end{document}